\documentclass[reqno, oneside]{amsart}
\usepackage{hyperref,url}
\usepackage{geometry}
\usepackage{pdfsync}
\usepackage{cite}

\theoremstyle{definition}

\newcommand\Atop[2]{\genfrac{}{}{0pt}{}{#1}{#2}}
\def\F(#1,#2,#3;#4){{}_2 F_1 \left ( \Atop{#1,#2} {#3} \, ; \, #4 \right )}
\def\Fthree(#1,#2,#3,#4,#5;#6){{}_3 F_2 \left ( \Atop{#1,#2,#3} {#4,#5} \, ; \, #6 \right )}
\def\f(#1,#2,#3;#4){F(#1,#2;#3;#4)}

\def\tfracx#1/#2@{\tfrac{#1}{#2}}

\usepackage[dvipsnames]{xcolor}

\begin{document}

\title{Zeilberger to the Rescue}

\author {Moa Apagodu}

\address{ Department of Mathematics, Virginia Commonwealth
University, Richmond, VA 23284}

\email{mapagodu@vcu.edu}

\begin{abstract} 
We provide both human and computer (even better collaboration between the two) proofs to four recent American Mathematical Monthly problems, namely problems \#11897,  \# 11899, \#11916, and  \#11928. We also show that problem 11928 may lead to interesting combinatorial identities. 
\end{abstract}

\maketitle

\begin{center}
Dedicated to Herbert S. Wilf, 1931-2012. 
\end{center}

\vskip .3in

\noindent We will demonstrate that Zeilberger's creative telescoping proof methods coupled with human touch proves most of monthly problems involving the binomial coefficients. Problem 11928 leads to the following curious identity involving the Catalan number.

$$
\sum_{i=0}^{n}\sum_{j=0}^{m}\binom{n}{i}\binom{m}{j}C(i+j) = \sum_{k=0}^{n+m}\binom{n+m}{k}C(k),
$$

\vskip .2in

\noindent where $C(n)$ is the $n$th Catalan number.
\vskip .3in

\noindent {\bf Problem \#11897}. {\it Proposed by P. Dalyay, Szeged, Hungary.\rm Prove for $n\geq0$, that}

$$
\sum_{\substack{k+l=n\\ k, \,l\geq0}} \frac1{k+1}\binom{2k}k\binom{2l+2}{l+1}=2\binom{2n+2}n.
$$

\vskip .2in

\noindent {\bf First solution using generating functions}: First we recall a theorem from product of power series, namely \\

\noindent {\it Theorem [S. H. Wilf,  Generatingfunctionology, p. 36] If $f$ and $g$ are ordinary power series generating functions for sequences $\{a_n\}$ and $\{b_n\}$, then $fg$ is the ordinary power series generating function for the sequence }

$$
\left \{\sum_{s+t=n; s, t \geq 0} a_sb_t\right\}_{n=0}^{\infty} \,\,.
$$

\noindent Our solution will make use of the following well known formulas (\cite{W}, pp 52-54)    

\begin{eqnarray}
\frac{1-\sqrt{1-4x}}{2x} = \sum_{k=0}^{\infty} \frac{1}{k+1}\binom{2k}k x^k
\end{eqnarray}

\begin{eqnarray}
 \frac{1}{\sqrt{1-4x}} = \sum_{k=0}^{\infty} \binom{2k}{k} x^k
\end{eqnarray}

\begin{eqnarray}
\frac{1}{\sqrt{1-4x}}\left(\frac{1-\sqrt{1-4x}}{2x}\right)^k &=& \sum_{n=0}^{\infty} \binom{2n+k}{n} x^n
\end{eqnarray}

\noindent From (2), with change of summation variable, we get 

\begin{eqnarray}
\sum_{j=0}^{\infty} \binom{2j+2}{j+1} x^j= \sum_{j=1}^{\infty} \binom{2j}{j} x^{j-1}=\frac{1}{x}\left(\frac{1}{\sqrt{1-4x}}-1\right) = \frac{1-\sqrt{1-4x}}{x\sqrt{1-4x}}\,\,.
\end{eqnarray}

\noindent Therefore, combining (1) and (4)  with the theorem, the right-side of the sum in question has ordinary power series generating function 

$$
2 \frac{1}{\sqrt{1-4x}}\left(\frac{1-\sqrt{1-4x}}{2x}\right)^2\,\,.
$$

\noindent Finally, the identity follows from (3).

\vskip .2in 

\noindent {\bf Second (pocket size)  proof using Gosper's decision procedure \cite{G}}.  This time rewrite the sum in the form

$$
\sum_{k=0}^{n}\frac1{k+1}\binom{2k}k\binom{2(n-k+1)}{n-k+1}=2\binom{2n+2}n.
$$

\noindent and let $F(n,k)$ be the summand. By Gosper's algorithm, the hypergeometric term  

$$
G(k)={\frac { \left( -2\,n+2\,k-3 \right)  \left( k+1 \right) {2\,k\choose 
k}{2\,n-2\,k+2\choose n-k+1}k!}{ \left( k+1 \right) !\, \left( n+2
 \right) }}
$$

\noindent is an anti-difference of $F(n,k)$, that is,

$$
F(n,k)=G(k+1)-G(k).
$$

\noindent Now adding both sides over $k$ for $0 \leq n \leq n$, we end up with the identity above, namely

$$
\sum_{k=0}^{n}\frac1{k+1}\binom{2k}k\binom{2(n-k+1)}{n-k+1}= G(n+1)-G(0)=2\binom{2n+2}n.
$$


\vskip .4in
\noindent  {\bf Problem \#11899}. {\it Proposed by J. Sorel, Romania. \rm Show that for any positive integer $n$},
$$\sum_{k=0}^n\binom{2n}k\binom{2n+1}k+\sum_{k=n+1}^{2n+1}\binom{2n}{k-1}\binom{2n+1}k=\binom{4n+1}{2n}+\binom{2n}n^2.$$

\vskip .2in

\noindent We start by observing that the second sum on the right-hand side is equal to the first sum. To see this, re-write the second sum as 

\begin{eqnarray}
\sum_{k=n+1}^{2n+1}\binom{2n}{k-1}\binom{2n+1}k = \sum_{k=n+1}^{2n+1}\binom{2n}{2n+1-k}\binom{2n+1}{2n+1-k}\,\,,
\end{eqnarray}

\noindent and make the change of variable $m = 2n+1-k$ to obtain

\begin{eqnarray}
\sum_{k=n+1}^{2n+1}\binom{2n}{k-1}\binom{2n+1}k = \sum_{m=0}^{n}\binom{2n}{m}\binom{2n+1}{m}\,\,.
\end{eqnarray}

\noindent Therefor, the identity to be shown is equivalent to 

\begin{eqnarray}
\sum_{k=0}^n2\binom{2n}k\binom{2n+1}{k}=\binom{4n+1}{2n}+\binom{2n}n^2.
\end{eqnarray}

\noindent If $w(n)$ is the sum on the left-hand side, then application of Zelbeger's creative telescoping method \cite{Z} (go to Maple and type $ZeilbergerRecurrence(F(n,k),n,k,w,0..n)$, where $F(n,k)$ is the summand on the left-hand side of (3) )  yields that the sum satisfies the non-homogeneous linear recurrence

$$
 (2n^2+5n+3)w(n+1)-(32n^2+64n+30)w(n)=-\frac{(n+1)(16n^2+38n+18)}{n^2} {2n \choose n+1}^2
.$$

\noindent Now it is a routine exercise to show that the right-hand side also satisfies this recurrence. Verify that both sides equal to 2 for $n=0$ to complete the proof. 

\vskip .4in

\noindent {\bf Problem \#11916}. {\it Proposed by Hideyuki Ohtsuka, Saitama, Japan, and Roberto Tauraso, Universita di Roma " Tor Vergata," Rome, Italy. Show that if $n$, $r$, and $s$ are positive integers, then }

$$
{ n+r \choose n} \sum_{k=0}^{s-1} { r+k \choose r -1}{ n +k\choose n} = { n+s \choose n} \sum_{k=0}^{r-1} { s+k \choose s -1}{ n+k\choose n} \,\,.
$$

\noindent  First re-write the identity as 

$$
\sum_{k=0}^{s-1} { n+r \choose n} { r+k \choose r -1}{ n +k\choose n} =\sum_{k=0}^{r-1}  { n+s \choose n}  { s+k \choose s -1}{ n+k\choose n} \eqno{(1)}\,.
$$

\noindent We use the Wilf-Zeilberger Method to prove the identity. Denote the sum on the left-side of (1) by $f(n)$ and on the right-side by $g(n)$. Then, $f(n)$ and $g(n)$ are solutions of the first-order non-homogeneous difference equation 

$$ nw(n)-(n+1)w(n+1) = - {n+s \choose n} {s+r\choose s} {n+r\choose r} \frac{sr}{n+1}\,\,.$$

\noindent  To see this, call the summand on the left-side of (1) $F1(n,k)$ and on the right-side $F2(n,k)$. Also define two companion functions 
$$
G1(n,k):=F1(n,k)\frac{(k+1)k}{n+1}\,,
$$

\noindent and 

$$
G2(n,k):=F2(n,k)\frac{(k+1)k}{n+1}\,\,.
$$

\noindent Then, first check that $nF1(n,k)-(n+1)F1(n+1,k)= G1(n,k+1)-G1(n,k)$ and $nF2(n,k)-(n+1)F2(n+1,k)= G2(n,k+1)-G2(n,k)$, and sum the first of these equations from $k=0$ to $k=s-1$ and the second equation from $k=0$ to $k=r-1$. Now show that $G1(n,0)=G2(n,0)=0$ and $G1(n,s) = G2(n,r)$, which equals the right-side of the non-homogeneous difference equation. This establishes that $f(n)$ and $g(n)$ satisfy the difference equation.\\

\noindent Finally, since both $f(n)$ and $g(n)$ satisfy the first-order difference equation with the initial condition $f(1)=g(1)= r(r+1){r+s \choose s-1}$, we must have $f(n)=g(n)$ for all $n \geq 1$ and any positive integers $r$ and $s$.

\vskip .4in

\noindent {\bf Problem \#11928}. {\it Proposed by Hideyuki Ohtsuka, Saitama, Japan. For positive integers $n$ and $m$ and for a sequence $<a_i>$, prove }

$$
\sum_{i=0}^{n}\sum_{j=0}^{m}\binom{n}{i}\binom{m}{j}a_{i+j} = \sum_{k=0}^{n+m}\binom{n+m}{k}a_k,
$$

\noindent and 

$$
\sum_{i<j}\binom{n}{i}\binom{n}{j}\binom{i+j}{n} = \sum_{i<j}\binom{n}{i}\binom{n}{j}^2.
$$

\noindent For the first identity,  using Vandermonde's convolution, we can rewrite the single sum on the right side as

\begin{eqnarray}
\sum_{k=0}^{n+m}\binom{n+m}{k}a_k &=&  \sum_{k=0}^{n+m}\sum_{i=0}^{k}\binom{n}{i}\binom{m}{k-i}a_k   \\
& = & \sum_{i=0}^{n+m}\sum_{k=i}^{m+n}\binom{n}{i}\binom{m}{k-i}a_k\\
&=&\sum_{i=0}^{n+m}\sum_{j+i=i}^{m+n}\binom{n}{i}\binom{m}{j}a_{i+j} \\
&=& \sum_{i=0}^{n+m}\sum_{j=0}^{n+m-i}\binom{n}{i}\binom{m}{j}a_{i+j}\\
&=& \sum_{i=0}^{n}\sum_{j=0}^{n+m-i}\binom{n}{i}\binom{m}{j}a_{i+j}\\
&=& \sum_{i=0}^{n}\sum_{j=0}^{m}\binom{n}{i}\binom{m}{j}a_{i+j}.
\end{eqnarray}

\noindent The second equality is by reversing the order of summation;  the third equality is by change of variable of summation ( $j=k-i$); and  (5) and (6) follow from $\binom{n}{k}=0$ for $k > n$. This completes the proof of the first identity.

\newpage

\noindent If we take $a_i=\binom{i}{n}$, then the first identity leads to  

$$
\sum_{j=0}^{n}\sum_{i=0}^{m}\binom{n}{i}\binom{m}{j}\binom{i+j}{n}=\sum_{k=0}^{n+m}\binom{n+m}{k}\binom{k}{n}.
$$

\noindent  Using the identity $\binom{n}{k}\binom{k}{m}= \binom{n}{m}\binom{n-k}{m-k}$, the right side evaluates to

\begin{eqnarray*}
\sum_{k=0}^{n+m}\binom{n+m}{k}\binom{k}{n}&=&\sum_{k=0}^{2n}\binom{n+m}{k}\binom{k}{n}\\
                                                               &=&\sum_{k=0}^{n+m}\binom{n+m}{n}\binom{n+m-n}{k-n}\\
                                                               &=&\sum_{k=0}^{n+m}\binom{n+m}{n}\binom{m}{k-n}\\
                                                               &=&\binom{n+m}{n}\sum_{j=0}^{m}\binom{m}{j}\\
                                                               &=&\binom{n+m}{n}2^m.
\end{eqnarray*}

\noindent  This gives the following nice identity: For positive integers $m$ and $n$, 

\begin{equation}
\sum_{j=0}^{n}\sum_{i=0}^{m}\binom{n}{i}\binom{m}{j}\binom{i+j}{n}=\binom{n+m}{n}2^m . 
\end{equation}

\vskip .2in

\noindent To prove the second identity, taking $m=n$ in (7), we get

\begin{equation}
\sum_{j=0}^{n}\sum_{i=0}^{n}\binom{n}{i}\binom{n}{j}\binom{i+j}{n}=\sum_{k=0}^{2n}\binom{2n}{k}\binom{k}{n}.
\end{equation}

\noindent The left side  of (8)  can be written as

$$
\sum_{j=0}^{n}\sum_{i=0}^{n}\binom{n}{i}\binom{n}{j}\binom{i+j}{n}= \sum_{0\leq i<j\leq n}\binom{n}{i}\binom{n}{j}\binom{i+j}{n}+
$$
$$
\sum_{0\leq j<i\leq n}\binom{n}{i}\binom{n}{j}\binom{i+j}{n}+\sum_{0\leq i=j\leq n}\binom{n}{i}\binom{n}{j}\binom{i+j}{n}.
$$

\noindent Using the symmetry in $i$ and $j$, we can simplify this sum to  

$$
\sum_{j=0}^{n}\sum_{i=0}^{n}\binom{n}{i}\binom{n}{j}\binom{i+j}{n}=2 \sum_{0\leq i<j\leq n}\binom{n}{i}\binom{n}{j}\binom{i+j}{n}+\sum_{i=0}^{n}\binom{n}{i}^2\binom{2i}{n}.
$$

\noindent On the other hand  using symmetry in $i$ and $j$ and the identity $\displaystyle{\binom{2n}{n}=\sum_{j=0}^{n} \binom{n}{j}^2}$,  we can write the right side of (8) as

\begin{eqnarray*}
\sum_{k=0}^{2n}\binom{2n}{k}\binom{k}{n} &=& \binom{2n}{n}2^n\\
                                                                          &=&\sum_{i=0}^{n}\sum_{j=0}^{n}\binom{n}{i}\binom{n}{j}^2\\
                                                                          &=&2\sum_{0\leq i<j\leq n}\binom{n}{i}\binom{n}{j}^2+\sum_{0\leq i=j\leq n}\binom{n}{i}\binom{n}{j}^2\\
                                                                          &=&2\sum_{0\leq j<i\leq n}\binom{n}{i}\binom{n}{j}^2+\sum_{i=0}^{n}\binom{n}{i}^3.
\end{eqnarray*}

\noindent To complete the proof of the second identity, we must shown that 

$$
\sum_{i=0}^{n}\binom{n}{i}^2\binom{2i}{n}=\sum_{i=0}^{n}\binom{n}{i}^3.
$$

\noindent We accomplish that we appeal to Zelbeger's creative telescoping method\cite{Z}. Denote the left and right side by $a_n$ and  $b_n$ respectively. Then both sequences start with $1, 2, 10, 56, 346, 2252$ and satisfy the second order recurrence (computed using the Zeilberger algorithm)

$$
(n+2)^2 w(n+2)-(7 n^2+21 n+16) w(n+1)-8( n+1^2)w(n)=0.
$$

\noindent Therefore,  $a(n)=b(n)$ for all positive integers $n$. This completes the poof of the second identity.

\vskip .2in

\noindent  Remark: This shows that with careful choice of $\{a_i\}$, one can obtain (perhaps a nontrivial) binomial identities. For example, if we take $a_i$ is the $i^{th}$ Catalan number, then we get 

$$
\sum_{i=0}^{n}\sum_{j=0}^{m}\binom{n}{i}\binom{m}{j}\frac{1}{i+j+1}{2(i+j) \choose i+j} = \sum_{k=0}^{n+m}\binom{n+m}{k}\frac{1}{k+1} { 2k \choose k}.
$$

\vskip .3in

\noindent {\bf Remark } Proofs of Problem 11899 and Problem 11916 are also provided in \cite{A} as a special case of a general theorem. Here we provided direct proof to these problems.

\end{document}